\documentclass[preprint,12pt]{elsarticle}

\usepackage{algorithmic}
\usepackage{amsfonts}
\usepackage{amsmath}
\usepackage{amssymb}
\usepackage{bbm}
\usepackage{color}
\usepackage{eurosym}
\usepackage{graphicx}
\usepackage{latexsym}
\usepackage[T1]{fontenc}
\usepackage{upgreek}
\usepackage{mathtools}

\newtheorem{theorem}{Theorem}
\newtheorem{corollary}{Corollary}
\newtheorem{lemma}{Lemma}

\newtheorem{example}{Example}

\newtheorem{remark}{Remark}

\newcommand{\bfm}{\mbox{$\mbox{\boldmath $m$}$}} 
\newcommand{\bfn}{\mbox{$\mbox{\boldmath $n$}$}} 

\newcommand{\sbfm}{\mbox{$\mbox{\boldmath \scriptsize $m$}$}} 
\newcommand{\sbfn}{\mbox{$\mbox{\boldmath \scriptsize $n$}$}} 

\newcommand{\bfzero}{\mbox{$\mbox{\boldmath $0$}$}} 
\newcommand{\bfone}{\mbox{$\mbox{\boldmath $1$}$}}

\newcommand{\bigzero}{\scalebox{0.9}{\mbox{\large $0$}}}
\newcommand{\sbfzero}{\scalebox{0.7}{\bfzero}}

\date{}

\journal{Applied Mathematics Letters}

\begin{document}

\begin{frontmatter}

\title{An explicit Maclaurin series solution to a classic non-autonomous abstract evolution equation}

\author[tas,ab]{Andrew Bassom}
\author[scg,ph]{Phil Howlett\corref{cor1}}
\author[mel,pt]{Peter Taylor}

\address[tas]{School of Natural Sciences, University of Tasmania, Hobart, Tasmania, 7001, Australia}
\address[scg]{Scheduling and Control Group, Centre for Industrial and Applied Mathematics, UniSA STEM, University of South Australia, Mawson Lakes 5095, Australia}
\address[mel]{School of Mathematics and Statistics, University of Melbourne, Parkville, Victoria, 3010, Australia}

\address[ab]{email:~andrew.bassom@utas.edu.au. url:~https://orcid.org/0000-0003-3275-7801}
\address[ph]{email:~phil.howlett@unisa.edu.au. url:~http://orcid.org/0000-0003-2382-8137}
\address[pt]{email:~taylorpg@unimelb.edu.au. url:~https://orcid.org/0000-0001-7600-5383}

\cortext[cor1] {Corresponding author}

\begin{abstract}
\label{abs}We propose and justify an explicit Maclaurin series solution to a classic non-autonomous abstract evolution equation for bounded linear operators on Banach space.
\end{abstract}


\end{frontmatter}

\section{Introduction}
\label{s:int}

Let $X$ be a Banach space and consider the abstract non-autonomous evolution equation
\begin{equation}
\label{naee}
dR(t)/dt = A(t)R(t)
\end{equation}
with $R(0) = I$ where the kernel $A \in {\mathcal C}^{\infty}([0,a), {\mathcal B}(X))$ is well defined if $t \in [0,a) \subseteq {\mathbb R}$ for some $a > 0$.  Suppose there exists a solution $R \in {\mathcal C}^{\infty}([0,r), {\mathcal B}(X))$ for some $r \in [0,a)$ such that $R(t)$ is real analytic for $t \in (-r, r)$.  Then we use basic power series manipulation to prove that $A(t) = A_0 + A_1t + \cdots$ is necessarily real analytic when $t \in (-b, b)$ for some $b >0$.  Next we apply the initial condition $R(0) = I$ and the recursive relationship 
\begin{equation}
\label{rnrr}
nR_n = A_0R_{n-1} + \cdots + A_{n-1}R_0
\end{equation}
for all $n \in {\mathbb N} = \{1,2,3,\ldots \}$ to show that $R(t) = I + R_1 t + R_2t^3 + \cdots$ is also analytic for $t \in (-b, b)$.  Finally we prove that  
\begin{equation}
\label{rna}
R_n = \mbox{$\sum_{q=0}^{n-1}$} \left[ \mbox{$\sum_{\sbfm \in S_{n,q}}$} \pi_{\sbfm} A_{\sbfm} \right] \in {\mathcal B}(X)
\end{equation}
for each $n \in {\mathbb N}$ where we define the index set
\begin{equation}
\label{snq}
S_{n,q} = \{ \bfm = (m_1,\ldots,m_{n-q}) \in ({\mathbb N}-1)^{n-q} = \{ 0,1,2,\ldots\}^{n-q} \mid \langle \bfone, \bfm \rangle = q\}
\end{equation}
for each $n \in {\mathbb N}$ and $q \in \{0,1,\ldots, n-1 \}$ and the coefficient
\begin{equation}
\label{pim} 
\pi_{\sbfm} = \mbox{$\prod_{j=1}^{n-q}$} \left[ n-q+1-j +\left( \mbox{$\sum_{s=j}^{n-q}$} m_s \right) \right]^{-1} \in {\mathbb R}
\end{equation}
for each $\bfm \in S_{n,q}$.  We write $A_{\sbfm} = A_{m_1} A_{m_2} \cdots A_{m_{n-q}}$ for each $\bfm \in S_{n,q}$ and $\bfone = \bfone_{n-q} = (1,\ldots,1) \in ({\mathbb N}-1)^{n-q}$.  It follows that $\langle \bfone, \bfm \rangle = q$ is the total index for the product $A_{\sbfm}$ and $n-q$ is the total exponent.  The explicit non-recursive formula (\ref{rna}) is our main result.  Various forms of the recursive relationship (\ref{rnrr}) are well known \cite{kuk1} but a non-recursive formula for $R_n$ when $n \in {\mathbb N}$ seems to be more elusive.  

\begin{example}
\label{ex1}

{\rm To illustrate (\ref{rna}) we calculate the coefficient $R_6$.  It is convenient to write the coefficient formula (\ref{pim}) in the equivalent more intuitive form
$$
\pi_{\sbfm} =  1/[n(n-1-m_1)(n-2-(m_1+m_2))\cdots (q+1 - (m_1+\cdots +m_{n-q-1}))]
$$
for each $\bfm = (m_1,m_2,\ldots,m_{n-q}) \in S_{n,q}$.  When $n=6$ we calculate 
\begin{eqnarray*} 
R_6 & = & \frac{A_0^6}{720} + \left( \frac{A_0^4A_1}{720} + \frac{A_0^3A_1A_0}{360} +\frac{A_0^2A_1A_0^2}{240} + \frac{A_0A_1A_0^3}{180} + \frac{A_1A_0^4}{144} \right) \\
& & + \left( \frac{A_0^3A_2}{360} + \frac{A_0^2A_1^2}{240} + \frac{A_0^2A_2A_0}{120} + \frac{A_0A_1A_0A_1}{180} + \frac{A_0A_1^2A_0}{90} \right. \\
& & \hspace{1.5cm} \left. + \frac{A_0A_2A_0^2}{60} + \frac{A_1A_0^2A_1}{144} + \frac{A_1A_0A_1A_0}{72} + \frac{A_1^2A_0^2}{48} + \frac{A_2A_0^3}{36} \right) \\
& & + \left( \frac{A_0^2A_3}{120} + \frac{A_0A_1A_2}{90} + \frac{A_0A_2A_1}{60} + \frac{A_0A_3A_0}{30} +    \frac{A_1A_0A_2}{72} \right. \\
& & \hspace{1.5cm} \left. + \frac{A_1^3}{48} + \frac{A_1A_2A_0}{24} + \frac{A_2A_0A_1}{36} + \frac{A_2A_1A_0}{18} + \frac{A_3A_0^2}{12} \right) \\
& & + \left( \frac{A_0A_4}{30} + \frac{A_1A_3}{24} + \frac{A_2^2}{18} + \frac{A_3A_1}{12} + \frac{A_4A_0}{6}  \right) + \frac{A_5}{6}
\end{eqnarray*}
where the relevant subsets are $S_{6,q}$ for each $q \in \{0,1,\ldots, 5 \}$.  For $S_{6,2}$ we have
\begin{itemize}
\item $\bfm = (0,0,0,2)$: $\pi_{(0,0,0,2)} = \frac{1}{6} \frac{1}{5} \frac{1}{4} \frac{1}{3} = \frac{1}{360}$.
\item \rule{0cm}{0.4cm}$\bfm = (0,0,1,1)$: $\pi_{(0,0,1,1)} = \frac{1}{6} \frac{1}{5} \frac{1}{4} \frac{1}{2} = \frac{1}{240}$.
\item \rule{0cm}{0.4cm}$\bfm = (0,0,2,0)$: $\pi_{(0,0,2,0)} = \frac{1}{6} \frac{1}{5} \frac{1}{4} \frac{1}{1} = \frac{1}{120}$.
\item \rule{0cm}{0.4cm}$\bfm = (0,1,0,1)$: $\pi_{(0,1,0,1)} = \frac{1}{6} \frac{1}{5} \frac{1}{3} \frac{1}{2} = \frac{1}{180}$.
\item \rule{0cm}{0.4cm}$\bfm = (0,1,1,0)$: $\pi_{(0,1,1,0)} = \frac{1}{6} \frac{1}{5} \frac{1}{3} \frac{1}{1} = \frac{1}{90}$.
\item \rule{0cm}{0.4cm}$\bfm = (0,2,0,0)$: $\pi_{(0,2,0,0)} = \frac{1}{6} \frac{1}{5} \frac{1}{2} \frac{1}{1} = \frac{1}{60}$.
\item \rule{0cm}{0.4cm}$\bfm = (1,0,0,1)$: $\pi_{(1,0,0,1)} = \frac{1}{6} \frac{1}{4} \frac{1}{3} \frac{1}{2} = \frac{1}{144}$.
\item \rule{0cm}{0.4cm}$\bfm = (1,0,1,0)$: $\pi_{(1,0,1,0)} = \frac{1}{6} \frac{1}{4} \frac{1}{3} \frac{1}{1} = \frac{1}{72}$.
\item \rule{0cm}{0.4cm}$\bfm = (1,1,0,0)$: $\pi_{(1,1,0,0)} = \frac{1}{6} \frac{1}{4} \frac{1}{2} \frac{1}{1} = \frac{1}{48}$.
\item \rule{0cm}{0.4cm}$\bfm = (2,0,0,0)$: $\pi_{(2,0,0,0)} = \frac{1}{6} \frac{1}{3} \frac{1}{2} \frac{1}{1} = \frac{1}{36}$.
\end{itemize}
The other calculations are similar.  The number of distinct indices $\bfm \in S_{n,q}$ is the number $\binom{n-1}{q}$ of non-negative solutions to the diophantine equation $m_1 + \cdots+ m_{n-q} = q$ for each $q \in \{0,1,\ldots,n-1\}$.  The total number of terms in $R_n$ is $2^{n-1} = \sum_{q=0}^{n-1} \binom{n-1}{q}$.}  $\hfill \Box$
\end{example}

\section{Literature review}
\label{s:lr}

Let $X$ be a complex Banach space and suppose that the evolution kernel ${\mathfrak A}(t) \in {\mathcal B}(X)$ is a bounded linear operator on $X$ for all $t \in [a,b]$.  Kato and Tanabe~\cite{kat1, tan1} and Sobolevski~\cite{sob1} found conditions on ${\mathfrak A}(t)$ that guarantee the existence of an evolutionary operator ${\mathfrak R}(s,t)$ such that 
${\mathfrak R}(t,s){\mathfrak R}(s,r) = {\mathfrak R}(t,r)$ for $s \in [r,t] \subseteq [a,b]$ with ${\mathfrak R}(t,t) = I$ and such that ${\mathfrak R}(t,s)$ satisfies $\partial {\mathfrak R}(t,s)/ \partial t = {\mathfrak A}(t) {\mathfrak R}(t,s)$ subject to the initial condition ${\mathfrak R}(s,s) = I$.  For an extended discussion see Yosida~\cite[Section XIV.5, pp 440\textendash 445]{yos1}.  Other more recent papers \cite{acq1,dig1,gis1,lat1} have continued this work.  

The Peano\textendash Baker series \cite{bak1,cam1,pea1} is a long-standing computational technique.  If $A(t) \in {\mathcal B}({\mathbb C}^k)$ is continuous for $t \in [a,b]$ and $\int_{[a,b]} \|A(t)\|dt < \infty$ the series is defined by $R(t,s) = \sum_{n \in {\mathbb N}-1} U_n(t,s)$ with $U_0(t,s) = I$ and $U_n(t,s) = \int_{[s,t]} A(\sigma) U_{n-1}(\sigma,s) d\sigma$ for all $[s,t] \subseteq [a,b]$ and all $n \in {\mathbb N}$.  In this case \cite[Theorem 2, p 158]{baa1} the series is compactly convergent on $[a,b]$ and $\partial R(t,s)/\partial t = A(t)R(t,s)$ with $R(s,s) = I$.  If $A(t) \in {\mathcal B}(X)$ is real analytic at $t=0$ the iteration for $R(t) = R(0,t)$ exactly reproduces the first $n+1$ terms of the Maclaurin series after $n$ iterations but does not reveal an explicit formula for the coefficients. 

Although not directly related to the methods used here the Magnus expansion \cite{bla1, mag1} seeks an infinite series  solution to (\ref{naee}) in the form of an operator exponential.  The expansion is especially useful for applications relating to classical and quantum mechanics because the truncated series often shares important qualitative properties with the exact solution.  The expansion is computed iteratively because explicit formul{\ae} are notoriously complicated \cite[Section 2.2, p 20]{bla1}.  The work by Magnus \cite{mag1} followed fundamental work on Lie algebras relating to the {\em so-called} CBHD exponential theorem of Campbell\textendash Baker\textendash Hausdorff\textendash Dynkin.  Achilles and Bonfiglio \cite{ach1} give a detailed history of this famous theorem.  An explicit formula for CBHD did not appear until the 1947 paper by Dynkin \cite{dyn1}.  

\section{Preliminaries}
\label{s:pre}

To prove our main result we assume that the kernel $A(t) \in {\mathcal B}(X)$ and the solution $R(t) \in {\mathcal B}(X)$ are both real analytic functions at $t = 0$.  We justify this assumption by arguing that $R(t)$ is real analytic at $t=0$ if and only if $A(t)$ is real analytic at $t = 0$.  We note two standard results to begin.

\begin{lemma}
\label{lem:1}
Let $F \in {\mathcal C}^{\infty}([0,a], {\mathcal B}(X))$ for some $a \in {\mathbb R}$ with $a > 0$.  The function $F(t)$ is real analytic for $t \in [0,b]$ where $b \in (0,a] \subseteq {\mathbb R}$ if and only if there exist $M, c > 0$ with $\sup_{t \in [0,b]}\| F^{(k)}(t)/k! \| \leq M / c^k$ for all $k \in {\mathbb N}-1$. $\hfill \Box$
\end{lemma}

{\bf Proof.} We simply observe that $F(t) = \sum_{k \in {\mathbb N}-1} F^{(k)}(0)t^k/k!$ for all $t \in (-b, b)$ if and only if $\| F^{(k)}(0) \| c^k/k!  \rightarrow 0$ as $k \rightarrow \infty$ for all $c \in (0, b)$.  See \cite{hik1,whi1} for further discussion. $\hfill \Box$

\begin{corollary}
\label{cor:1}
If $F(t)$ is real analytic for $t \in (-b,b)$ for some $b > 0$ the extended function defined by $F(z) \coloneqq \sum_{k \in {\mathbb N}-1} F^{(k)}(0)z^k/k!$ is analytic for all $z \in {\mathbb C}$ with $\lvert z \rvert < b$. $\hfill \Box$
\end{corollary}

The next two results are specific to our particular problem.

\begin{lemma}
\label{lem:2}
Let $A(t), R(t) \in {\mathcal B}(X)$ satisfy {\rm (\ref{naee})}.  If the solution $R(t)$ is real analytic at $t=0$ then the kernel $A(t)$ is real analytic at $t=0$.  $\hfill \Box$
\end{lemma}

{\bf Proof.} There exists $r > 0$ such that $R(t) = I + \sum_{k \in {\mathbb N}} R^{(k)}(0)t^k/k!$ for all $t \in (-r,r)$.  For each $s \in (0,r)$ there is some $M = M_s > 0$ such that $\| R^{(k)}(0) \| \leq M_s k!/s^k$ for all $k \in {\mathbb N}$.  Therefore the extended function $R(z) \coloneqq I + \sum_{k \in {\mathbb N}} R^{(k)}(0)z^k/k!$ is analytic for $z \in {\mathbb C}$ with $\lvert z \rvert < r$.  Write $R(z) = I + K(z)$ where $K(z) \coloneqq R_1z + R_2z^2 + \cdots$ for $\lvert z \rvert < r$ and $\|K(z)\| \rightarrow 0$ as $\lvert z \rvert \rightarrow 0$.  Hence there is some $b > 0$ such that
$$
[I + K(z)]^{-1} = I - K(z) + K(z)^2 - \cdots = I - R_1z + (R_1^2-R_2)z^2 + \cdots
$$
is well defined and analytic for $z \in {\mathbb C}$ with $\lvert z \rvert < b$.  Therefore (\ref{naee}) implies
$$
A(t) = [R_1 + 2R_2t + 3R_3t^2 + \cdots][I - R_1t + (R_1^2-R_2)t^2 + \cdots] = A_0 + A_1t + A_2t^2 + \cdots
$$
for all $t \in (-b, b)$.  It follows that $A(t)$ is real analytic at $t = 0$.  $\hfill \Box$

\begin{lemma}
\label{lem:3}
Let $\{A_j\}_{j \in {\mathbb N}-1} \in {\mathcal B}(X)$ and suppose that $A(t) \coloneqq A_0 +A_1t + A_2t^2 + \cdots$ is real analytic on some nontrivial interval $t \in (-b, b)$.  Let $R_0 = I \in {\mathcal B}(X)$ and use the recursive relationship
\begin{equation}
\label{rnrr+}
nR_n = A_0R_{n-1} + \cdots + A_{n-1}R_0
\end{equation} 
for each $n \in {\mathbb N}$ to define a countable collection of coefficients $\{R_j\}_{j \in {\mathbb N}-1} \in {\mathcal B}(X)$.  The function $R(t) \coloneqq R_0 + R_1t + \cdots$ is real analytic  and satisfies the non-autonomous abstract evolution equation {\rm (\ref{naee})} for all $t \in (-b, b)$.  If $A(t)$ is a polynomial this result extends to all such finite intervals.  Thus it  is valid for all $t \in{\mathbb R}$.  $\hfill \Box$
\end{lemma}

{\bf Proof.} If we define $\alpha_n = M/b^n$ for each $n \in {\mathbb N}-1$ and $\{\rho_n\}_{n \in {\mathbb N}-1}$ by setting $\rho_0 = 1$ and $n \rho_n =  \alpha_0 \rho_{n-1} + \cdots +  \alpha_{n-1}\rho_0$ for each $n \in {\mathbb N}$ then $\rho_n = Mb (Mb + 1) \cdots (Mb + n-1)/(n! \, b^n)$.  This identity is easily proved by induction.  We claim that $\| R_n \| \leq  \rho_n$ for each $n \in {\mathbb N}$.  The claim is true for $n=1$ because $\| R_1 \| = \| A_0 \| = M = \rho_1$.  Suppose the claim is true for $n \leq m-1$.  It follows from (\ref{rnrr+}) that
$$
\| m R_m \| \leq \| A_0\| \| R_{m -1} \| + \cdots + \| A_{m-1} \| \| R_0\| \leq \alpha_0 \rho_{m-1} + \cdots +  \alpha_{m-1}\rho_0 = m \rho_m
$$
which means the claim is true for $n \leq m$.  Thus the desired result is true for all $n \in {\mathbb N}$.  Hence $R(t) \coloneqq R_0 + R_1t + \cdots$ is real analytic for $t \in (-b, b)$ with $\| R(t) \| \leq M/(1 - \lvert t \rvert/b)^{Mb}$.  A straightforward argument now shows that $R(t)$ satisfies (\ref{naee}) for all $t \in (-b, b)$.  $\hfill \Box$

\subsection{Some additional notation}
\label{s:not}

If $A(t)$ is a polynomial of degree $p \in {\mathbb N}$ then for each $q \in \{0,1,\ldots,\lfloor pn/(p+1) \rfloor \}$ where $k = \lfloor x \rfloor$ denotes the largest integer $k \leq x$ then the restricted index set of order $p$ is 
\begin{equation}
\label{snqp}
S_{n,q,p} = \{ \bfm \in \{0,1,\ldots,p\}^{n-q} \mid \langle \bfone, \bfm \rangle = q\} \subseteq S_{n,q}.
\end{equation}

\section{Statement of the main results}
\label{s:smr}

Let $X$ be a complex Banach space and let ${\mathcal B}(X)$ denote the space of bounded linear operators on $X$.  The following results will be proved in Section \ref{s:jmr}.

\begin{theorem}
\label{t1:smr}
Let $\{A_j\}_{j\in {\mathbb N}-1} \in {\mathcal B}(X)$ be an infinite collection of bounded linear operators with $\limsup_{j \rightarrow \infty} \|A_j\|^{1/j} = 1/b$ for some $b \in (0,\infty) \cup \{\infty\}$ and let $A(t) = \sum_{j \in {\mathbb N}-1} A_jt^j \in {\mathcal B}(X)$ for all $t \in (-b,b)$ be the corresponding analytic function.   For each $n \in {\mathbb N}$ and $q \in \{0,1,\ldots,n-1\}$ let $S_{n,q}$ be the index set defined by {\rm (\ref{snq})} and for each $\bfm \in S_{n,q}$ let $\pi_{\sbfm}$ be the coefficient defined by {\rm (\ref{pim})}.  The analytic function $R(t) = I + \sum_{j \in {\mathbb N}} R_jt^j$ where
\begin{equation}
\label{rnanal}
R_n = \mbox{$\sum_{q=0}^{n-1}$} \left[ \mbox{$\sum_{\sbfm \in S_{n,q}}$} \pi_{\sbfm} A_{\sbfm} \right] \in {\mathcal B}(X)
\end{equation}
for each $n \in {\mathbb N}$ is the unique solution to {\rm (\ref{naee})} for all $t \in (-b,b)$.   $\hfill \Box$
\end{theorem} 

\begin{theorem}
\label{t2:smr}
Let $\{A_j\}_{j=1}^p \in {\mathcal B}(X)$ be a finite collection of bounded linear operators and let $A(t) = \sum_{j=0}^p A_jt^j \in {\mathcal B}(X)$ for all $t \in {\mathbb R}$ be the corresponding polynomial.  For each $n \in {\mathbb N}$ and $q \in \{0,1,\ldots,\lfloor pn/(p+1) \rfloor \}$ let $S_{n,q,p}$ be the restricted index set of order $p$ defined by {\rm (\ref{snqp})} and for each $\bfm \in S_{n,q,p}$ let $\pi_{\sbfm}$ be the coefficient defined by {\rm (\ref{pim})}. The analytic function $R(t) = I + \sum_{j \in {\mathbb N}} R_jt^j$ where 
\begin{equation}
\label{rnpoly}
R_n = \mbox{$\sum_{q=0}^{ \lfloor pn/(p+1) \rfloor }$} \left[ \mbox{$\sum_{\sbfm \in S_{n,q,p}}$} \pi_{\sbfm} A_{\sbfm} \right] \in {\mathcal B}(X)
\end{equation}
for all $n \in {\mathbb N}$  is the unique solution to {\rm (\ref{naee})} for all $t \in {\mathbb R}$.  Note that $\lfloor pn/(p+1) \rfloor = n-k$ when $n \in [k(p+1), (k+1)(p+1))$ for all $k \in {\mathbb N}$.  $\hfill \Box$
\end{theorem}

\section{Justification of the main results}
\label{s:jmr}

We may assume from Lemma \ref{lem:3} that $A(t) = A_0 + A_1t + \cdots$ and $R(t) = I + R_1t + R_2t^2 + \cdots$ are analytic and that $R(t)$ is a solution to (\ref{naee}) for all $t \in (-b,b)$. 

\begin{lemma}
\label{lem4}
For each $n \in {\mathbb N}$ the coefficient $R_n$ is given by the formula
\begin{equation}
\label{bsrn}
R_n = \mbox{$\sum_{q=0}^{n-1}$} \left[ \mbox{$\sum_{\sbfm \in S_{n,q}}$} \pi_{\sbfm} A_{\sbfm} \right]
\end{equation}
where $S_{n,q} = \{ \bfm \in ({\mathbb N}-1)^{n-q} \mid \langle \bfone, \bfm \rangle = q \}$ for each $q \in \{0,1,\ldots,n-1\}$ and the constant $\pi_{\sbfm}$ is defined by the formula
\begin{equation}
\label{bspim}
\pi_{\sbfm} = \mbox{$\prod_{j=1}^{n-q}$} \left[ n-q+1-j +\left( \mbox{$\sum_{s=j}^{n-q}$} m_s \right) \right]^{-1}
\end{equation}
for all $\bfm \in S_{n,q}$.  In {\rm(\ref{bsrn})} we have written $A_{\sbfm} = A_{m_1} \cdots A_{m_{n-q}}$. $\hfill \Box$
\end{lemma}

{\bf Proof.}  When $n=1$ we have $S_{1,0} = \{ \bfm = \bfzero = (0) \}$.  Therefore $m_1 = 0$ and so $\pi_{\sbfzero} = 1/(m_1+1) = 1$.  Since $R_0 = I$ the proposed formula gives $R_1 = A_0$.  This agrees with the value determined by the recursive relationship (\ref{rnrr+}).  Therefore (\ref{bsrn}) and (\ref{bspim}) are true for $n = 1$.  Choose $k \in {\mathbb N} +1$ and suppose (\ref{bsrn}) and (\ref{bspim}) are true for $n < k$.  For each $j \in \{0,1,\ldots,k-2\}$ the inductive hypothesis gives
$$
R_{k-1-j} = \mbox{$\sum_{q=0}^{k-2-j}$} \left( \mbox{$\sum_{\sbfn \in S_{k-1-j, q}}$} \pi_{\sbfn} A_{\sbfn} \right) = \mbox{$\sum_{p=j}^{k-2}$} \left( \mbox{$\sum_{\sbfn \in S_{k-1-j, p-j}}$} \pi_{\sbfn} A_{\sbfn} \right)
$$
because $\bfn = (n_1,\ldots,n_{k-1-q-j}) \in S_{k-1-j, q} \iff \bfn = (n_1,\ldots,n_{k-1-p}) \in S_{k-1-j, p-j}$ if $p = q + j$.  The recursive relationship $kR_k = A_0R_{k-1} + \cdots + A_{k-1}R_0$ can be rewritten as
$$
kR_k = A_{k-1} + \mbox{$\sum_{j=0}^{k-2}$} A_j \mbox{$\sum_{p=j}^{k-2}$} \left( \mbox{$\sum_{\sbfn \in S_{k-1-j, p-j}}$} \pi_{\sbfn} A_{\sbfn} \right).
$$
If we reverse the order of summation for the two outermost sums we get
$$
kR_k = A_{k-1} + \mbox{$\sum_{p=0}^{k-2}$} \left\{ \mbox{$\sum_{j=0}^{k-2}$} \left( \mbox{$\sum_{\sbfn \in S_{k-1-j, p-j}}$} \pi_{\sbfn} A_jA_{\sbfn} \right) \right\}.
$$
Because $\bfn \in S_{k-1-p, p-j} \iff (j,\bfn) \in S_{k,p}$ for each $j =0,\ldots,k-2$ it follows that
$$
R_k = A_{k-1}/k + \mbox{$\sum_{p=0}^{k-2}$}  \left( \mbox{$\sum_{(j,\sbfn) \in S_{k,p}}$} \pi_{\sbfn} A_{(j,\sbfn)}/k \right)
$$
where we note that $\pi_{(k-1)} = 1/k$ and $\pi_{\sbfm} = \pi_{(j,\sbfn)} = \pi_{\sbfn}/k$ for each $\bfm = (j,\bfn) \in S_{k,p}$.  Now the inductive hypothesis with $n = k-1$ and $q = p$ in (\ref{bspim}) shows that
$$
\pi_{\sbfm} = (1/k)  \mbox{$\prod_{r=1}^{k-1-p}$} \left[  k-p -r +\left(\,  \mbox{$\sum_{s=r}^{k-1-p}$} n_s \right) \right]^{-1}.
$$
If we substitute $t = r+1$ and use the fact that $n_s = m_{s+1}$ for $s=1,\ldots,k-1-p$ we get
$$
\pi_{\sbfm} = \mbox{$\prod_{t=1}^{k-p}$} \left[  k-p+1-t +\left(\,  \mbox{$\sum_{s=t}^{k-p}$} m_s \right) \right]^{-1}.
$$
Therefore (\ref{bsrn}) and (\ref{bspim}) are true for $\bfm \in S_{k,p}$ and all $ p \in \{0,1,\ldots,k-1\}$.  It follows by induction that these formul{\ae} are also true when $\bfm \in S_{n,q}$ for all $n \in {\mathbb N}$ and $q \in \{0,1,\ldots,n-1\}$.   $\hfill \Box$ 

It is well known \cite[Theorem 2, p 158]{baa1} that the evolution equation has a unique solution under much weaker conditions than those assumed here.  Consequently our main result, Theorem~\ref{t1:smr}, follows from Lemma \ref{lem4}.  However, in our case, it is also true that term by term differentiation of the explicit formula for $R(t)$ shows that $R(t)$ satisfies the evolution equation and hence shows that a solution exists.  Of course the explicit formula {\em per se} shows that the solution is unique.

\vspace{0.25cm}
We now turn our attention to Theorem \ref{t2:smr} where $A(t)$ is a polynomial.  If $A_j = \bigzero$ for $j \in {\mathbb N}+p = \{p+1,p+2,p+3,\ldots\}$ then for each $n \in {\mathbb N}$ we wish to find the largest value of $q$ for which the coefficient
$R_{n,q} = \sum_{\sbfm \in S_{n,q}} \pi_{\sbfm}A_{\sbfm}$ is nonzero.   Every nonzero term that contributes to $R_{n,q}$ is identified by an index $\bfm \in S_{n,q}$ with $m_1 + \cdots + m_{n-q} = q$ and $m_j \in \{0,1,\ldots,p\}$ for all $j \in \{1,\ldots,n-q\}$.  The largest value of $q$, or equivalently the smallest possible value of $n - q$ that allows a nonzero contribution, will occur when we have the largest possible number of indices $\{m_j\}_{j=1}^{n-q}$ equal to $p$.  If $\bfm = p \bfone_{n-q}$ then
\begin{equation}
\label{prod1}
\langle \bfone_{n-q}, \bfm \rangle = q  \Longrightarrow p(n - q) = q
\end{equation}
If $\bfm = (p \bfone_{\ell},r,p\bfone_{n-q-1-\ell})$ for some $\ell \in \{0,\ldots,n-q-1\}$ and $r \in \{1,\ldots,p-1\}$ then
\begin{equation}
\label{prod2}
\langle \bfone_{n-q}, \bfm \rangle = q \Longrightarrow p(n-q-1)+r = q.
\end{equation}
We can cover both cases at once if we also allow $r = p$ in (\ref{prod2}).  Now it follows that
$$
q =  \left(p(n-1) + r \right)/ (p+1) \Longrightarrow pn/(p+1) - 1 <  q \leq pn/(p+1).
$$
Because $q = q_{\max}$ is an integer we must have $q_{\max} = \lfloor pn/(p+1) \rfloor$.  This justifies defining a restricted index set $S_{n,q,p}$ in (\ref{snqp}).  For all $n \in {\mathbb N}$ and $q \in \{0,1,\ldots, \lfloor pn/(p+1) \rfloor\}$ we can replace the general formula (\ref{bsrn}) by the more discriminating formula
\begin{equation}
\label{bsrnp}
R_n  = \mbox{$\sum_{q=0}^{\lfloor pn/(p+1) \rfloor}$} \left[ \mbox{$\sum_{\sbfm \in S_{n,q,p}}$} \pi_{\sbfm} A_{\sbfm} \right].
\end{equation}
It is now a matter of setting $A_j = \bigzero$ for all $j >p$ in Theorem \ref{t1:smr} to justify Theorem~\ref{t2:smr}.  It is still true with $A(t) = A_0 + A_1t + \cdots + A_pt^p$ that we could have $A_j = \bigzero$ for some $j \in \{0,1,\ldots,p-1\}$.  Thus there may still be some zero products in the reduced expression $(\ref{bsrnp})$ and some non-zero coefficients $\pi_{\sbfm}$ for indices $\bfm \in S_{n,q,p}$ where $A_{\sbfm} = \bigzero$.

\begin{remark}\
\label{rem3}
We could eliminate all zero products by defining a list ${\mathcal J} = \{j_1,j_2,\ldots \} \subseteq {\mathbb N}-1$ of all nonzero coefficients $A_j$ and a restricted index set
$$
S_{n,q,{\mathcal J}} = \{ \bfm = (m_1,\ldots,m_{n-q}) \in {\mathcal J}^{n-q} \mid m_1 + \cdots + m_{n-q} = q \} \subseteq S_{n,q}
$$
for each $n \in {\mathbb N}$ and each $q \in \{0,1,\ldots,n-1\}$.  We will not pursue a detailed analysis.  $\hfill \Box$
\end{remark}

\section{Possible extensions}
\label{s:pe}
If the coefficients $\{A_n\}_{n \in {\mathbb N}-1}$ are known then the coefficients $\{R_n\}_{n \in {\mathbb N}-1}$ are uniquely determined by the initial condition and a linear recurrence relation.  This leads to an explicit non-recursive formula for the coefficients $\{R_n\}_{n \in {\mathbb N}}$.  If the coefficients $\{R_n\}_{n \in {\mathbb N}-1}$ are known then the coefficients $\{A_n\}_{n \in {\mathbb N}}$ are uniquely determined by a rearranged version of the same linear recurrence relation.  This suggests the possibility of an explicit non-recursive formula for the coefficients $\{A_n\}_{n \in {\mathbb N}}$.

\section{Conclusions}
\label{s:con}
We have shown that $A(t) = \sum_{n \in {\mathbb N}-1}A_nt^n$ is real analytic at $t = 0$ if and only if $R(t) = \sum_{n \in {\mathbb N}-1} R_nt^n$ is real analytic at $t = 0$ and we have found an explicit formula for $R_n$ in terms of $A_0, A_1,\ldots,A_{n-1}$ for all $n \in {\mathbb N}$.  In this case it is necessary and sufficient for convergence of the respective Maclaurin series when $t \in (-r,r)$ that there are constants $c, r > 0$ such that $\|A_n\|, \|R_n\| \leq c/r^n$ for all $n \in {\mathbb N}$.  It follows that $A(t)$ and $R(t)$ are also bounded and that the conditions on the kernel and solution are nicely balanced.   This is not true in general.   Consider the situation where $A(t)$ is smooth but not analytic.  If $A(t)$ is unbounded and densely defined on some Hilbert space $H$ and skew adjoint near $t=0$ and if $R(0)=I$ an elementary argument shows that
$$
(dR^*/dt) R = R^*A^*R = - R^*AR = - R^*(dR/dt) \Rightarrow d(R^*R)/dt = \bigzero \Rightarrow R^*R = I.
$$
Thus $R(t)$ is unitary.  It follows that $R(t)$ is bounded.  However $R(t)$ is not analytic because, if it were, our earlier arguments would show that $A(t)$ must also be analytic and bounded.  Having said that it may be interesting to investigate the case where the bounded linear operator $A(t) \in {\mathcal B}(H)$ is real analytic and skew adjoint near $t=0$. 

\section{Acknowledgements}
\label{s:ack}
We would like to thank an anonymous referee who alerted us to the fact that $R(t)$ is unitary when $A(t)$ is skew adjoint.

\section{Declaration}
\label{s:dec}

The authors declare that they have no conflicting interests relating to this paper.

\end{document}